\documentclass[11pt]{amsart}
\usepackage{geometry}                
\usepackage{graphicx}
\usepackage{amssymb}
\RequirePackage{ifthen}
\usepackage{amsmath}
\usepackage{amsthm}
\usepackage{xcolor}
\usepackage{epstopdf}
\newtheorem{theorem}{Theorem}

\newcommand{\enmom}{{m}}
\newcommand{\benmom}{{\mathbf{m}}}

\newenvironment{talk}[4][]
{\noindent\parbox{\textwidth}{%
\begin{center}%
\textbf{\boldmath #3}\\[\smallskipamount]\textsc{#2}%
{\\%
(joint work with #1)}%
\end{center}}\par\nopagebreak\bigskip\nopagebreak}%
{\bigskip\bigskip\goodbreak}

\begin{document}


\begin{talk}[Erwann Delay \cite{ChDelayHPETv1} and Gregory Galloway~\cite{ChGallowayHPET}]{Piotr T. Chru\'sciel}
{Hyperbolic positive energy theorems}
{Chru\'sciel, Piotr T.}

\noindent

It is convenient to start this report with a few definitions. We say that a Riemannian manifold $(M,g)$ is \emph{conformally compact} if there exists a compact manifold with boundary $\widehat  M$ such that the following holds:
First, we allow $M$ to have a boundary, which is then necessarily compact. Next,
$M$ is the interior of $\widehat  M$, whose boundary is the union of the boundary of $M$ and of a number of new boundary components, at least one, which form the \emph{conformal boundary at infinity}. Further, there exists on $\widehat  M$ a smooth function
$\Omega\ge 0 $  which is positive on $M$, and which vanishes precisely on the new boundary components of $\widehat  M$,
with $d\Omega$ nowhere vanishing there. Finally,   the tensor field
$\Omega^2 g$ extends to a smooth metric on $\widehat  M$.

We will say that a conformally compact manifold $(M,g)$ is  \emph{asymptotically locally hyperbolic} (ALH) if all sectional curvatures approach minus one as the conformal boundary  at infinity is approached.
An ALH metric has an \emph{asymptotically hyperbolic} (AH) component of its boundary at infinity, or an AH end, if the conformal metric on that component of the boundary is conformal to a round sphere.

A useful global invariant of an ALH-but-not-AH end is its mass $m$, while for AH ends we have an energy-momentum vector $\benmom \equiv (\enmom_\mu)$
 \cite{ChHerzlich,ChNagyATMP,Wang}
 (compare~\cite{AbbottDeser,ChruscielSimon}).
For this one considers metrics $g$ which asymptote, at a suitable rate, to a background metric $\mathring g$. It is assumed that $\mathring g$ admits nontrivial \emph{static potentials} which,
in dimension $n$,
are defined as solutions of the overdetermined system of equations
 \begin{equation}\label{12VIII1711-}
    \mathring D_i \mathring D_ j V = \big(  \mathring R_{ij} -   \frac { \mathring R}{n-1}  \mathring g_{ij} \big) V
       \,.
 \end{equation}
Here $\mathring D$ is the covariant derivative of the background metric $\mathring g$, while $\mathring R_{ij}$ is its Ricci tensor and $\mathring R $ is the trace $\mathring g^{ij} \mathring R_{ij}$.
To every static potential $V$ and asymptotic end $\partial M$ one associates a mass $m=m(V,\partial M)$  by the formula
  \cite{HerzlichRicciMass} (compare \cite[Equation~(IV.40)]{BCHKK})
\begin{eqnarray}
 m(V,\partial M)
    &= &
   - 
   \lim_{x\rightarrow0}\int_{\{x\} \times \partial M}   D^j V
    ( R{}^i{}_j - \frac {R{}}{n}\delta^i_j)
    \,
    d\sigma_i
     \,,
           \label{18IX20.4}
    \end{eqnarray}
where
$R_{ij}$ is the Ricci tensor of the metric $g$, $R$ its trace, and
we have ignored an overall dimension-dependent positive multiplicative factor which is often used in the physics literature. Here $\partial M $ is a component of conformal infinity, and $x$ is a coordinate near $\partial M$ so that $\partial M$ is given by the equation $\{x=0\}$.

The difference between AH ends and general ALH ends arises from the dimension of the space of static potentials. Indeed, the AH case is the only one where this dimension is larger than one. Then   $\mathring g$ is taken to  be the hyperbolic metric, which can be  written as the following metric on $\mathbb{R}^n$:
\begin{equation}\label{3XI18.8}
\mathring g =
   \frac {dr^2 }{1+r^2} + r^2 {d\Omega^2_{n-1}}
 \,,
\end{equation}
where ${d\Omega^2_{n-1}}$ is the unit round metric on $S^{n-1}$. In this coordinate system
 a basis of the space of static potentials is provided by the functions
$$
 V_{0} = \sqrt{r^2+1}
 \,,
 \qquad
 V_{i}=    x^i
  \,.
$$
One defines the components $m_\mu$ of the \emph{energy-momentum vector} $\benmom$ as
\begin{equation}\label{15VII21.1}
  m_\mu := m(V_{\mu})
  \,.
\end{equation}
One checks that $\benmom$ transforms as a Lorentz vector under conformal transformations of $S^{n-1}$, so that its Lorentzian norm is a geometric invariant.

For all remaining ALH ends, the mass is directly an invariant~\cite{ChNagyATMP}.

Strictly speaking, a rescaling of $V$ by a constant is always possible, and a preferred scale can be set as follows: In AH ends the standard normalisation is the one just described.
In all remaining cases, in any chosen ALH end we can write $\mathring g$ as
\begin{equation}\label{16VII21.2}
  \mathring g = x^{-2}(dx^2 + \mathring h)
  \,,
  \quad
  \mathring h(\partial_x, \cdot) = 0
  \,,
\end{equation}
with the volume of $\partial M$, calculated in the metric $\mathring h|_{x=0}$, normalised to one. One then normalises $V$ so that $\lim_{x\to 0}xV=1$.

There is a closely related definition of energy-momentum for  \emph{asymptotically flat} general relativistic initial data sets $(M,g,K)$ which, perhaps somewhat unexpectedly, turns out to be relevant for the \emph{asymptotically hyperbolic} problem at hand, and which is invoked in Theorem~\ref{T12XII18.1} below, we refer the reader to~\cite{Bartnik86,ChErice} for details.

In the case of spherical conformal infinity,  it has been known that $\benmom$ is timelike future pointing under a spin condition~\cite{ChHerzlich,CJL,GHHP,Maerten,Wang}, or under restrictive hypotheses~\cite{AnderssonGallowayCai,CGNP}.
 In my talk in Oberwolfach, summarised here, I  reported on results presented in \cite{ChDelayHPETv1,ChGallowayHPET} where it is shown how to remove these hypotheses.

The starting point of the analysis in \cite{ChDelayHPETv1} is the following result:

\begin{theorem}
  \label{T12XII18.1}
  Let $({M},g)$ be an asymptotically Euclidean Riemannian manifold, where ${M}$ is the union of a compact set and of an asymptotically flat region, of dimension $n\ge 3$. Suppose that the general relativistic initial data set $({M},g,K)$ possesses a well defined energy-momentum vector $\benmom$. If the dominant energy condition holds, then
 $\benmom$ is timelike
    future pointing or vanishes.
  Furthermore,  in the last case $({M},g,K)$ arises from a hypersurface in Minkowski spacetime.
\end{theorem}

A published proof of Theorem~\ref{T12XII18.1} in dimensions less than or equal to seven can be found in~\cite{Eichmair:PET,EHLS,HuangLee}, building upon~\cite{SchoenYauPNAS,SchoenYauPMT1,SchoenYauPMT2}. A proof covering all dimensions is available in preprint form in \cite{Lohkamp2}, with the borderline cases covered in \cite{ChBeig1,ChMaerten,HuangLee}.
Conjecturally, this result also follows in all dimensions basing on the preprint \cite{SchoenYau2017}.

In~\cite{ChDelayHPETv1} it is shown how
Theorem~\ref{T12XII18.1}, together with   the perturbation results in~\cite{CGNP} and the gluing constructions of~\cite{ChDelayExotic}, can be used to remove all unnatural restrictions in the proof of positivity of asymptotically hyperbolic mass:

\begin{theorem}\label{rigid}
Let
$({M},g)$ denote an $n$-dimensional Riemannian manifold  which is the union of a compact set and an AH end. If the
 scalar curvature $R(g)$
satisfies
$
 R(g)\ge -n(n-1)
  \,,
$
then the energy-momentum vector of $(M,g)$ is causal future pointing, or vanishes.
\end{theorem}

The impossibility of a \emph{null} future pointing energy-momentum vector, under the hypotheses above, has been established in~\cite{HuangJangMartin}.

Theorem~\ref{rigid} has been generalised in~\cite{ChGallowayHPET} to allow manifolds with several ends, and with boundaries satisfying an optimal mean-curvature condition:

\begin{theorem}
 \label{T26V21.1}
Let $(M,g)$ be a conformally compact
  $n$-dimensional, $3\le n \le 7$,  asymptotically locally hyperbolic   manifold  with boundary.
Assume that the scalar curvature  of $M$ satisfies $R(g) \ge  -n(n-1)$, and that the boundary has mean curvature
 $H \le n-1$
 with respect to the normal pointing into $M$.  Then,   the energy-momentum vector
 ${\bf m}$ of every spherical component of the conformal boundary at infinity of $(M,g)$ is future causal.
\end{theorem}

In this theorem neither the boundary $\partial M$, nor the conformal boundary at infinity of $M$, need to be connected.
The proof relies heavily on the results of~\cite{Eichmair:2020wbj}, which assume $3\le n \le 7$.

The above theorems concern AH ends, and one is led to wonder about properties of mass for ALH-but-not-AH ends. Here the following is known: First, positivity is known on manifolds with suitable spin structure~\cite{Wang}, or under restrictive
conditions~\cite{ChGallowayHPET,CGNP}, but such $(M,g)$ are scarce.  Next, boundaryless conformally compact examples with negative mass and toroidal infinity are due to Horowitz and Myers~\cite{HorowitzMyers}; nontrivial quotients of spheres at infinity with, again, negative mass have been constructed by Chen and Zhang~\cite{ChenZhang}. Finally, a natural negative lower bound, together with a Penrose-type inequality (compare~\cite{ChruscielSimon,GibbonsPenrose}), has been established by Lee and Neves in~\cite{LeeNeves} for a class of three dimensional models with higher genus conformal infinity.

\bibliographystyle{amsplain}

\begin{thebibliography}{10}

\bibitem{AbbottDeser}
L.F. Abbott and S.~Deser, \emph{Stability of gravity with a cosmological
  constant}, Nucl.\ Phys. \textbf{B195} (1982), 76--96.

\bibitem{AnderssonGallowayCai}
L.~Andersson, M.~Cai, and G.J. Galloway, \emph{Rigidity and positivity of mass
  for asymptotically hyperbolic manifolds}, Ann.\ H.~Poincar\'e \textbf{9}
  (2008), 1--33, arXiv:math.dg/0703259. \MR{MR2389888 (2009e:53054)}

\bibitem{Bartnik86}
R.~Bartnik, \emph{The mass of an asymptotically flat manifold}, Commun.\ Pure
  Appl.\ Math. \textbf{39} (1986), 661--693. \MR{849427 (88b:58144)}

\bibitem{BCHKK}
H.~Barzegar, P.T. Chru\'{s}ciel, and M.~{H\"orzinger}, \emph{Energy in
  higher-dimensional spacetimes}, Phys.\ Rev.\ D \textbf{96} (2017), 124002, 25
  pp., arXiv:1708.03122 [gr-qc].

\bibitem{ChBeig1}
R.~Beig and P.T. Chru\'{s}ciel, \emph{Killing vectors in asymptotically flat
  spacetimes: {I. A}symptotically translational {K}illing vectors and the rigid
  positive energy theorem}, Jour.\ Math.\ Phys. \textbf{37} (1996), 1939--1961,
  arXiv:gr-qc/9510015.

\bibitem{ChenZhang}
J.~Chen and X.~Zhang, \emph{Metrics of {E}guchi-{H}anson types with the
  negative constant scalar curvature}, Jour.\ Geom.\ Phys. \textbf{161} (2021),
  104010, 10, arXiv:2007.15964 [math.DG]. \MR{4180104}

\bibitem{ChErice}
P.T. Chru\'{s}ciel, \emph{Boundary conditions at spatial infinity from a
  {H}amiltonian point of view}, Topological Properties and Global Structure of
  Space--Time (P.\ Bergmann and V.\ de~Sabbata, eds.), Plenum Press, New York,
  1986, pp. 49--59, arXiv:1312.0254 [gr-qc].

\bibitem{ChDelayExotic}
P.T. Chru\'{s}ciel and E.~Delay, \emph{Exotic hyperbolic gluings}, Jour.\
  Diff.\ Geom. \textbf{108} (2018), 243--293, arXiv:1511.07858 [gr-qc].

\bibitem{ChDelayHPETv1}
P.T. Chru\'{s}ciel and E.~Delay, \emph{{The hyperbolic positive energy
  theorem}},  (2019), arXiv:1901.05263v1 [math.DG].

\bibitem{ChGallowayHPET}
P.T. Chru\'sciel and G.J. Galloway, \emph{{Positive mass theorems for
  asymptotically hyperbolic Riemannian manifolds with boundary}}, Class. Quantum
  Grav. \textbf{38} (2021), 237001, arXiv:2107.05603 [gr-qc].

 

\bibitem{CGNP}
P.T. Chru\'{s}ciel, G.J. Galloway, L.~Nguyen, and T.-T. Paetz, \emph{{On the
  mass aspect function and positive energy theorems for asymptotically
  hyperbolic manifolds}}, Class. Quantum Grav. \textbf{35} (2018), 115015,
  arXiv:1801.03442 [gr-qc].

\bibitem{ChHerzlich}
P.T. Chru\'{s}ciel and M.~Herzlich, \emph{The mass of asymptotically hyperbolic
  {R}iemannian manifolds}, Pacific Jour.\ Math. \textbf{212} (2003), 231--264,
  arXiv:math/0110035 [math.DG]. \MR{MR2038048 (2005d:53052)}

\bibitem{CJL}
P.T. Chru\'{s}ciel, J.~Jezierski, and S.~Leski, \emph{The{ Trautman-Bondi} mass
  of hyperboloidal initial data sets}, Adv.\ Theor.\ Math.\ Phys. \textbf{8}
  (2004), 83--139, arXiv:gr-qc/0307109. \MR{MR2086675 (2005j:83027)}

\bibitem{ChMaerten}
P.T. Chru\'{s}ciel and D.~Maerten, \emph{Killing vectors in asymptotically flat
  spacetimes: {II}. {A}symptotically translational {K}illing vectors and the
  rigid positive energy theorem in higher dimensions}, Jour.\ Math.\ Phys.
  \textbf{47} (2006), 022502, arXiv:gr-qc/0512042. \MR{MR2208148 (2007b:83054)}

\bibitem{ChNagyATMP}
P.T. Chru\'{s}ciel and G.~Nagy, \emph{The mass of spacelike hypersurfaces in
  asymptotically {anti -- de Sitter} spacetimes}, Adv.\ Theor.\ Math.\ Phys.
  \textbf{5} (2001), 697--754, arXiv:gr-qc/0110014.

\bibitem{ChruscielSimon}
P.T. Chru\'{s}ciel and W.~Simon, \emph{Towards the classification of static
  vacuum spacetimes with negative cosmological constant}, Jour.\ Math.\ Phys.
  \textbf{42} (2001), 1779--1817, arXiv:gr-qc/0004032.

\bibitem{Eichmair:PET}
M.~Eichmair, \emph{{The Jang equation reduction of the spacetime positive
  energy theorem in dimensions less than eight}}, Commun.\ Math.\ Phys.
  \textbf{319} (2013), 575--593, arXiv:1206.2553 [math.dg]. \MR{3040369}

\bibitem{Eichmair:2020wbj}
M.~Eichmair, G.J. Galloway, and A.~Mendes, \emph{{Initial data rigidity
  results}},  (2020), arXiv:2009.09527 [gr-qc].

\bibitem{EHLS}
M.~Eichmair, L.-H. Huang, D.A. Lee, and R.~Schoen, \emph{The spacetime positive
  mass theorem in dimensions less than eight}, Jour.\ Eur.\ Math.\ Soc.\ (JEMS)
  \textbf{18} (2016), 83--121, arXiv:1110.2087 [math.DG]. \MR{3438380}

\bibitem{GibbonsPenrose}
G.W. Gibbons, \emph{Some comments on gravitational entropy and the inverse mean
  curvature flow}, Class.\ Quantum Grav. \textbf{16} (1999), 1677--1687,
  arXiv:hep-th/9809167. \MR{1697098}

\bibitem{GHHP}
G.W. Gibbons, S.W. Hawking, G.T. Horowitz, and M.J. Perry, \emph{Positive mass
  theorem for black holes}, Commun.\ Math.\ Phys. \textbf{88} (1983), 295--308.

\bibitem{HerzlichRicciMass}
M.~Herzlich, \emph{Computing asymptotic invariants with the {R}icci tensor on
  asymptotically flat and asymptotically hyperbolic manifolds}, Ann.\ Henri
  Poincar\'e \textbf{17} (2016), 3605--3617, arXiv:1503.00508 [math.DG].
  \MR{3568027}

\bibitem{HorowitzMyers}
G.T. Horowitz and R.C. Myers, \emph{The {AdS/CFT} correspondence and a new
  positive energy conjecture for general relativity}, Phys.\ Rev.\ D
  \textbf{59} (1998), 026005, arXiv:hep-th/9808079.

\bibitem{HuangJangMartin}
L.-H. Huang, H.C. Jang, and D.~Martin, \emph{{Mass rigidity for hyperbolic
  manifolds}}, Commun.\ Math.\ Phys. (2019), 1--21, arXiv:1904.12010 [math.DG].

\bibitem{HuangLee}
L.-H. Huang and D.A. Lee, \emph{Equality in the spacetime {P}ositive {M}ass
  {T}heorem}, Commun.\ Math.\ Phys. \textbf{376} (2020), 2379--2407,
  arXiv:1706.03732 [math.DG]. \MR{4104553}

\bibitem{LeeNeves}
D.A. Lee and A.~Neves, \emph{The {P}enrose inequality for {a}symptotically
  {l}ocally {h}yperbolic spaces with nonpositive mass}, Commun.\ Math.\ Phys.
  \textbf{339} (2015), 327--352. \MR{3370607}

\bibitem{Lohkamp2}
J.~{Lohkamp}, \emph{{The Higher Dimensional Positive Mass Theorem II}},
  (2016), arXiv:1612.07505 [math.DG].

\bibitem{Maerten}
D.~Maerten, \emph{Positive energy-momentum theorem in asymptotically {anti-de
  Sitter} spacetimes}, Ann.\ H.Poincar\'e \textbf{7} (2006), 975--1011,
  arXiv:math.DG/0506061. \MR{MR2254757 (2007d:83016)}

\bibitem{SchoenYauPNAS}
R.~Schoen and S.-T. Yau, \emph{Complete manifolds with nonnegative scalar
  curvature and the positive action conjecture in general relativity}, Proc.\
  Nat.\ Acad\. Sci.\ U.S.A. \textbf{76} (1979), 1024--1025. \MR{MR524327
  (80k:58034)}

\bibitem{SchoenYauPMT1}
\bysame, \emph{On the proof of the positive mass conjecture in general
  relativity}, Commun.\ Math.\ Phys. \textbf{65} (1979), 45--76. \MR{MR526976
  (80j:83024)}

\bibitem{SchoenYauPMT2}
\bysame, \emph{Proof of the positive mass theorem. {II}}, Commun.\ Math.\ Phys.
  \textbf{79} (1981), 231--260. \MR{612249}

\bibitem{SchoenYau2017}
\bysame, \emph{{Positive Scalar Curvature and Minimal Hypersurface
  Singularities}},  (2017), arXiv:1704.05490 [math.DG].

\bibitem{Wang}
X.~Wang, \emph{Mass for asymptotically hyperbolic manifolds}, Jour.\ Diff.\
  Geom. \textbf{57} (2001), 273--299. \MR{MR1879228 (2003c:53044)}

\end{thebibliography}
\providecommand{\bysame}{\leavevmode\hbox to3em{\hrulefill}\thinspace}
\providecommand{\MR}{\relax\ifhmode\unskip\space\fi MR }
\providecommand{\MRhref}[2]{%
  \href{http://www.ams.org/mathscinet-getitem?mr=#1}{#2}
}
\providecommand{\href}[2]{#2}

\end{talk}
\end{document}